\documentclass[12pt,draft]{amsart}

\usepackage[english]{babel}
\usepackage{amsmath,amsthm,enumerate}
\usepackage{amsfonts,bm}
\newtheorem{thm}{Theorem}[section]
\newtheorem{cor}[thm]{Corollary}
\newtheorem{lem}[thm]{Lemma}

\newtheorem{defn}{Definition}[section]

\theoremstyle{definition}
\theoremstyle{remark}
\newtheorem{rem}[defn]{Remark}
\numberwithin{equation}{section}
\newtheorem{example}[defn]{Example}

\def\d{\partial}
\def\beq{\begin{equation}}
\def\eeq{\end{equation}}
\def\be{\begin{example}}
\def\ee{\end{example}}
\def\br{\begin{rem}}
\def\er{\end{rem}}
\def\bt{\begin{thm}}
\def\et{\end{thm}}
\def\L{\Lambda}
\def\l{\lambda}

\def\a{\alpha}
\def\b{\beta}
\def\z{{\bar z}}

\def\1#1{\overline{#1}}
\def\2#1{\widetilde{#1}}
\def\3#1{\widehat{#1}}
\def\4#1{\mathbb{#1}}
\def\5#1{\frak{#1}}
\def\6#1{{\mathcal{#1}}}

\def\bl{\begin{lem}}
\def\el{\end{lem}}

\def\bpf{\begin{proof}}
\def\epf{\end{proof}}
\def\ben{\begin{enumerate}}
\def\een{\end{enumerate}}
\def\beq{\begin{equation}}
\def\eeq{\end{equation}}

\newcommand{\R}{\mathbb{R}}
\newcommand{\Z}{\mathbb{Z}}
\newcommand{\Q}{\mathbb{Q}}
\newcommand{\C}{\mathbb{C}}
\newcommand{\N}{\mathbb{N}}
\newcommand{\B}{\ensuremath{\mathfrak B}}
\newcommand{\mt}{\ensuremath{\mathfrak M}}
\newcommand{\ct}{\ensuremath{\mathfrak C}}

\newcommand{\ls}{\ensuremath{\mathcal L}}
\DeclareMathOperator{\im}{Im}
\DeclareMathOperator{\re}{Re}

\date{\today}%

\begin{document}

\title[Sums of squares and torsion]{
	Sums of squares in pseudoconvex hypersurfaces
	and torsion phenomena for Catlin's boundary systems
}

\author[A. Basyrov]{Alexander Basyrov}

\address{Department of Mathematics, Statistics and Computer Science, University of Wisconsin-Stout, Menomonie, WI 54751}

\email{abasyrov@gmail.com}

\author[A. Nicoara]{Andreea C. Nicoara}

\author[D. Zaitsev]{Dmitri Zaitsev}

\address{School of Mathematics, Trinity College Dublin, Dublin 2, Ireland}

\email{anicoara@maths.tcd.ie, zaitsev@maths.tcd.ie}




\begin{abstract}
Given a pseudoconvex hypersurface in $\C^n$ and an arbitrary weight, 
we show the existence of local coordinates in which the polynomial model 
contains a particularly simple sum of squares of monomials.
Our second main result provides a normalization of a part of any Catlin boundary system.
We illustrate by an example that 
this normalization cannot be extended to the rest of the boundary system
due to the existence of what we refer to as torsion.
\end{abstract}

\maketitle

\section{Introduction}\label{intro}

The goal of this paper is to gain new geometric insight into the tools
developed for establishing global regularity and subelliptic estimates for the $\bar\d$-Neumann problem.
We are focusing on the approach by Catlin
in 
\cite{C84, C87}.
However, we expect our methods to 
also shed light
on Kohn's multiplier ideal technique
initiated in \cite{kohnacta}
and continued more recently in 
\cite{fribourgcatda, siunote,  ra, kimzaitsev, siunew}
as well as on
other research related to the $\bar\d$-equation
(see \cite[1.3]{dimatensors}) and potentially
 more general PDE's,
as evidenced by the program pioneered by Siu in \cite{siunew}.

Recall that Catlin established global regularity and subelliptic estimates
for the $\bar\d$-Neumann problem
as consequences of his Property (P) type conditions.
The only known proofs of Property (P) type conditions for
general smooth pseudoconvex finite type domains in the sense of D'Angelo~\cite{opendangelo}
rely on the techniques of multitype, polynomial models, and boundary systems
introduced in \cite{catlinbdry}.

One of the motivations in this paper is reducing the complexity in Catlin's techniques
by a more explicit use of pseudoconvexity and a precise normalization of the geometry.
%
%
Our first result is showing the existence of so-called positive balanced terms in 
polynomial models of pseudoconvex hypersurfaces.
We call a monomial $z_1^{\a_1} \z_1^{\b_1} \ldots z_n^{\a_n} \z_n^{\b_n} $
{\em balanced} if $\a_j=\b_j$ for all $j$.
We shall use the (standard) {\em lexicographic order} for multiweights  and
the {\em reverse lexicographic order} for multidegrees.
We prove the following:

\begin{thm}
  \label{main}
  Let $M\subset \C^n$ be a pseudoconvex real
  smooth hypersurface with $0\in M$.
  Let $z=(z_1, \ldots, z_n)$ be local 
  holomorphic coordinates in a neighborhood of $0$
  and
  $$
  	\mu= (\mu_1, \mu_2, \mu_3, \ldots, \mu_n),
	\quad
	1=\mu_1> \mu_2\ge  \mu_3 \ge \ldots \ge \mu_n > 0,
  $$  
   be a (multi-)weight such that $M$ is given by
$$
	r(z,\bar z) =0,
	\quad
	d r \ne 0,
	\quad
	r = O_\mu(1),
$$
i.e.\ the expansion of $r$  contains only terms of weight greater or equal to $1$.

Then after a weighted homogenous polynomial  change of coordinates,
the defining function 
$r$ of $M$ admits a decomposition
%
  \begin{equation}\label{main-eq}
    r(z,\bar z) =-2\re z_1 
  + p (z_2, \ldots, z_n, \bar z_2, \ldots, \bar z_n)
+ o_\mu(1),
  \end{equation}
  where 
  $p$ is a weighted homogeneous polynomial of weight $1$ 
containing, as part of its expansion, the sum of squares
\beq\label{bal}
    A_2 |z_2|^{2k_{22}} + A_3 |z_2|^{2k_{32}} |z_3|^{2k_{33}} 
    + \ldots 
    + A_n  |z_2|^{2k_{n2}} \cdots |z_n|^{2k_{nn}}
\eeq
with $A_j \ge 0$, $k_{jj} >0$ for all $j=2,\ldots, n,$
such that the (total) degree of $p$ in each $(z_j, \bar z_j)$ is not greater than $2k_{jj}$, 
and  $o_\mu(1)$ stands for a smooth function in 
$(\im z_1, z_2, \ldots, z_n, \bar z_2, \ldots, \bar z_n)$
whose formal Taylor series expansion
contains only terms
of weight greater than $1$.
In addition, for each $j$, the multidegree of each term of \eqref{bal} in 
$(z_2,\z_2), \ldots, (z_j,\z_j)$ is maximal among all
balanced monomials in $p$
in the reverse lexicographic order.

Furthermore, either all $A_j$ can be chosen positive, or the weight $\mu$ can be lowered
lexicographically.
\end{thm}

The following example illustrates that pseudoconvexity is an essential assumption:

\be
Let $M$ be the non-pseudoconvex domain given by $r<0$ with
$$
	r= -2\re z_1 + 2\re(z_2^2\z_3^3). 
$$
Then no weighted polynomial change of variables
transforms $M$ into a form satisfying the conclusion of Theorem~\ref{main}.
In fact, 
no terms as in \eqref{bal}
can be obtained (even with $A_j$ negative).
\ee

On the other hand, when $M$ is pseudoconvex,
a natural question arises whether
the statement could be  improved
 by reducing the polynomial $p$ to just a sum of 
 the terms in \eqref{bal}.
 It is not possible in general, however,
 as the following simple example illustrates: 
 
 \be
 Let $M\subset\C^2$ be the tube given by $r=0$ with
 $$
 	r=-2\re z_1 + (\re z_2)^2.
 $$
 Then the polynomial $(\re z_2)^2$
 is invariant under weighted homogeneous polynomial coordinate changes,
 and hence cannot be reduced to contain the terms \eqref{bal} only.
 \ee

 Furthermore, it is generally not possible to reduce $r$
 to the sum \eqref{bal} only,
 even when $r$ can be written as a sum of squares of holomorphic functions:
 
 \be
 Let $p,q \in \N,$ $p,q \ge 2.$ Consider 
\[
	r_0 =
	 -2\re z_1 + |z_2|^{2p}+|z_3|^{2q} + 
	 2\epsilon \re  z_2^p \bar z_3^q
\]
with $|\epsilon|<1$.
Then $r_0$ determines a real-algebraic pseudoconvex hypersurface 
that can be written as a sum of squares of holomorphic functions:
\beq\label{sq}
	r_0 = 
	-2\re z_1 +\left |z_2^p
	+{\epsilon}z_3^q\right|^2
	+ (1-\epsilon^2) | z_3|^{2q}.
\eeq
By a direct computation (or using e.g.\ \cite[Theorem~4.1]{martinIMRN}), 
it can be seen, however, that no
biholomorphic change of variables can transform $r_0$ into 
a sum of squares of the form 
$$
	|a_2z_2|^{2k_{22}} + A_3 |z_2|^{2k_{32}} |z_3|^{2k_{33}}.
$$
\ee

 The following example illustrates that 
 the terms \eqref{bal}
 cannot be in general reduced
 to sums of single powers $|z_j|^{2k_{jj}}$.
  \be
 Let $M\subset\C^3$ be given by $r=0$ with
 \beq\label{eqq}
 	r=-2\re z_1 + |z_2|^8 + |z_2|^4|z_3|^6,
 \eeq
a weighed homogeneous polynomial in $(z_1, z_2, z_3)$ and conjugates 
with corresponding weights $(1,8,12)$.
Then the only weighted homogeneous changes of coordinates
are the linear dilations $(z_1,z_2,z_3)\mapsto (a_1z_1,a_2z_2,a_3z_3)$
that clearly preserve \eqref{eqq} up to a change of coefficients.
In particular, it is not possible to obtain a conclusion similar to Theorem~\ref{main}
with \eqref{bal} consisting of sums of single powers $|z_j|^{2k_{jj}}$ only.
 \ee

Applying Theorem~\ref{main} to the Catlin multitype yields the following:

\begin{cor}
  \label{catlincase}
  Let $M$ be a pseudoconvex smooth real hypersurface in $\C^n$ with $0 \in M$ and 
  of the Catlin multitype 
  $$
  \Lambda= (1, \lambda_2, \lambda_3, \ldots, \lambda_n),
  \quad 
  \lambda_n < + \infty
  $$ 
  at $0$.
  
  Then there exists a holomorphic change of coordinates at $0$ 
  preserving the multitype so that the defining function for $M$ in the new coordinates is given by
  \begin{equation}\label{m-eq}
    r = -2\re z_1 + p(z_2, \ldots, z_n, \bar z_2, \ldots, \bar z_n) 
    + o_{\L^{-1}}(1),
  \end{equation}
    where 
    $p$ is a weighted homogeneous polynomial of weight $1$ that 
  contains the sum of squares
  \begin{equation}\label{a-terms}
  |z_2|^{2k_{22}} + |z_2|^{2k_{32}} |z_3|^{2k_{33}} + \ldots + |z_2|^{2k_{n2}} \cdots |z_n|^{2k_{nn}},
  \end{equation}
  where 
  $k_{jj} >0$ for all $j$ and the (total) degree of $p$ in each $(z_j, \bar z_j)$ is not greater than $k_{jj}.$ 


\end{cor}

\br
It follows from the weighted homogeneity of 
\eqref{a-terms} that
  \beq\label{lj}
   \sum_{l=2}^j \frac{2 k_{jl}}{ \lambda_l} =1,
   \quad
   j=1,\ldots, n.
  \eeq
  In particular, since $k_{jj}\ne 0$, the weights $\l_l$
  are uniquely determined from \eqref{lj}.
\er

A better understanding of the Catlin multitype
has been obtained when $M$ is a boundary
of a {\em convex} domain,
see \cite{bs, mcneal, jiyeyu}. 
However, even in this case,
Thereom~\ref{main}
seems to be new.


For a pseudoconvex hypersurface, the multitype and the commutator multitype coincide as Catlin established in \cite{catlinbdry}.  Without pseudoconvexity, however,
it might not be possible to obtain terms of the type 
$$
	|z_2|^{2k_{j2}} \cdots |z_n|^{2k_{jj}}
$$ 
via a change of variables, and the two multitypes may differ, as the following example due to Bloom shows:

 \be
 Consider 
\[
	r_0 =
	 \re z_1 + (\re z_2 +|z_3|^2)^2.
\]
The multitype at $0$ is given by $\mt=(1,2,4),$ whereas the commutator multitype at $0$ is lexicographically strictly larger: $\ct=(1,2, + \infty).$ Catlin proved that $\mt \leq \ct$ for any domain, and this example shows nothing more can be expected to hold in the absence of pseudoconvexity.
 \ee
 
Our second main result 
is a normalization of any boundary system for a pseudoconvex domain via a change of variables:

\bt
  \label{flatnessprop}
  Let $M$ be a pseudoconvex smooth real hypersurface in $\C^n$ with $0 \in M,$ Levi rank $s_0$ at $0,$ and 
  of the Catlin multitype 
  $$
  \Lambda= (1, \underbrace{2, \dots, 2}_{s_0}, \lambda_{s_0+2}, \ldots, \lambda_n)
  $$ 
  at $0$, where 
  $$
  	2<\lambda_{s_0+2}=\dots = \lambda_{s_0+s_1+1}< \lambda_{s_0+s_1+2} < +\infty.
$$ 
Then for any boundary system at $0$,
   $$
   	\B_n (0) = \{r_1, r_{s_0+2},\dots, r_n; L_2, \dots, L_n\},
$$ 
there exists a holomorphic change of coordinates at $0$ 
 preserving the multitype
 and  transforming $\B_n(0)$  into  
 $$
 	\tilde \B_n (0) = \{\tilde r_1,\tilde  r_{s_0+2}, \dots, \tilde r_n; \tilde L_2, \dots, \tilde L_n\},
$$ satisfying the normalization
  \begin{equation}\label{flatness-eq}
  \begin{aligned}
\tilde r_j &= \re z_j 
	+o(
		\lambda_j^{-1}
	), 
	\quad s_0+2 \leq j \leq s_0+s_1+1,\\
\tilde L_k&=\d_{z_k}
	+o(
		\lambda_k^{-1}
	), 
	\quad 2\le k\le s_0+s_1+2, \\
   \end{aligned}
  \end{equation}
   where the partial derivatives $\d_{\z_j}$ are counted with weight $-\l_j^{-1}$.
  \et

Therefore, the first function in the boundary system beside the defining function can always be brought to the simplest possible form.
In case several entries of the Catlin multitype are equal to the first $\lambda_j >2,$ 
all of their corresponding functions and vector fields in the boundary system can be normalized. This normalization process cannot be carried out on the subsequent functions, however. 
In fact, we provide a {\em counterexample} in the final section of the paper. 
This non-existence of a complete normalization provides a very important insight 
into the behavior of Catlin's boundary systems
that we call a {\em torsion phenomenon}.

In particular, for $M$ of Levi rank $0$ in $\C^3,$ 
this flattening result yields a simplified geometric picture, 
which we state as a corollary:

\begin{cor}
 \label{C3case}
   Let $M$ be a pseudoconvex smooth real hypersurface in $\C^3$ with $0 \in M,$ Levi rank $0,$ and 
  of the Catlin multitype 
  $$
  \Lambda= (1, \lambda_2, \lambda_3),
  \quad 
  \lambda_3 < + \infty
  $$ 
  at $0$. After a holomorphic change of variables, any boundary system at $0$  
  $$
  	\B_3 (0) = \{r_1, r_2, r_3; L_2, L_3\}
$$ becomes $$\tilde \B_3 (0) = \{\tilde r_1, \tilde r_2, \tilde r_3; \tilde L_2, \tilde L_3\}$$ with 
   \begin{equation}\label{flatness-eqc3}
  \begin{aligned}
\tilde r_{2} &= \re z_2 +o (\lambda_2^{-1}),\\
\tilde L_j&=\partial_{z_j}+
	o( \lambda_j^{-1} ), 
	\quad j =2,3,\\
   \end{aligned}
  \end{equation}
 \end{cor}

The paper is organized as follows: Section~\ref{examples} provides additional examples
illustrating various phenomena.
Section~\ref{notation} gives the relevant definitions and notation. Section~\ref{bdrysystems} defines Catlin's multitype and boundary systems. Section~\ref{estimates} presents some elementary auxiliary results such as one-dimensional estimates for non-negative homogeneous polynomials and a several variables version proven via scaling using the Newton polygon. Section~\ref{normalforms} contains the proofs of Theorem~\ref{main} and Corollary~\ref{catlincase} carried out in a sequence of lemmas. In the same section, it is shown that the first function in the boundary system can always be normalized, thus establishing Theorem~\ref{flatnessprop} and Corollary~\ref{C3case}. 
Section~\ref{r3counterexample} demonstrates by example that the same type of normalization cannot be carried out on subsequent functions in the boundary system.

\section{Further motivation and examples}\label{examples}

\be
The grandfather of all examples is a strongly pseudoconvex
hypersurface given by $r=0$ with
$$
	r=-2\re z_1 +  |z_2|^2 +\ldots + |z_n|^2 + o_{\L^{-1}}(1),
	\quad
	\L=(1, 2, \ldots, 2),
$$
where the whole leading polynomial (the quadric),
is already of the form \eqref{bal},
once diagonalized.
Of course, the diagonalization here is a special property of quadrics
that does not extend to higher degree polynomials.
\ee

More generally, 
it was shown by the third author in \cite{dimatensors}
that any pseudoconvex hypersurface admits the form
$r=0$ with
\beq\label{p4}
	r = -2\re z_1 + |z_2|^2 + \ldots + |z_{q-1}|^2 
	+ p_4(z_{[q+2,n]}, \z_{[q+2,n]} ) + o_{\L^{-1}}(1), 
\eeq
with the inverse weights 
$
	\L= (1, 2, \ldots, 2, 4, \ldots, 4),
$
where the number of $2$'s equals the Levi form rank $q$,
and $p_4$ represents the CR invariant quartic tensor
defined in \cite{dimatensors}.
In view of this general fact,
it suffices to only study terms
arising from $p_4$ and $o_{\L^{-1}}(1)$.

\be
Any real hypersurface $M\subset\C^2$ of finite type $m$
admits the form $r=0$ with
$$
	r = -2\re z_1 + p_m(z_2,\z_2) + o_{\L^{-1}}(1),
	\quad
	p_m = \sum_{j+k=m} a_{jk} z_2^j \z_2^{k},
$$
where $\L=(1,m)$ and $p_m$ is not harmonic.
In this (well-known) case, 
if $M$ is pseudoconvex,
Theorem~\ref{main}
implies that  $m=2l$ is even
and $p_m$ contains a nontrivial term $a_{ll} |z_2|^{2l}$,
$a_l\ne 0$.
Since $p_m$ is a (tensor) invariant,
it is clear that other terms $a_{jk} z_2^j \z_2^{k}$
with $j\ne k$ cannot be eliminated.
\ee

\be
For a general pseudoconvex hypersurface $M\subset\C^3$
of the form \eqref{main-eq}
in its multitype coordinates,
Corollary~\ref{catlincase} implies
that after a linear change of coordinates, 
$p$ contains nontrivial terms
\beq\label{p22}
	p = |z_2|^{2k_{22}} + |z_2|^{2k_{32}} |z_3|^{2k_{33}} + 
	\2p(z_2,z_3, \z_2,\z_3),
	\quad
	k_{22}, k_{33} \ge 1.
\eeq
where $\2p$ consists of all remaining terms.
In particular, the multitype $\L=(1, m_2, m_3)$
is uniquely determined from the identities
$$
	m_2= 2k_{22},
	\quad
	\frac{2k_{32}}{m_2} + \frac{2k_{33}}{m_3} = 1
$$
expressing the property that the first two terms in
\eqref{p22} are of weight $1$.

Furthermore, the additional degree property in Theorem~\ref{main}
asserts that the degree of $p$ 
in $(z_3,\z_3)$ equals $2k_{33}$. This property puts additional restrictions on the terms in $p$
and makes the choice of the first two terms in \eqref{p22} canonical.
For instance, if
$$
	p= |z_2|^4 + |z_2z_3|^2 + |z_3|^4,
$$
the degree condition forces us to choose the terms
$|z_2|^4 + |z_3|^4$ rather than $|z_2|^4 + |z_2z_3|^2$,
because the latter choice would violate
the property that the degree of $p$ in $(z_3,\z_3)$ is $2$.
\ee

\be
Let $M\subset\C^4$ be given by $r=0$ with
$$
	r = - 2\re z_1 + |z_2|^4 + |z_2|^2 |z_3|^2 + ( |z_2|^2 + |z_3|^2 ) |z_4|^2.
$$
Then the degree property in Theorem~\ref{main}
implies that the sum of squares in  \eqref{a-terms}
must be
\beq\label{sum0}
	|z_2|^4 + |z_2|^2 |z_3|^2 + |z_3|^2 |z_4|^2,
\eeq
since the remaining square term  
$|z_2|^2 |z_4|^2$
has the same degree in $(z_4,\z_4)$ as $|z_3|^2 |z_4|^2$
but lower degree in $(z_3,\z_3)$. 
That last restriction again determines
uniquely the terms in \eqref{sum0}.
\ee

\be
Corollary~\ref{catlincase} allows us to estimate how many inverse weights can arise as multitypes $\mt=(m_1, \dots,m_n)$ at $0$ of a pseudoconvex hypersurface $M\subset\C^n$ of finite $1$-type $m.$ We first recall that by part (4) of Catlin's Main Theorem in \cite{catlinbdry}, $m_n \leq m.$ It is always the case that $m_1=1,$ $m_2$ is an integer, and $2 \leq m_2 \leq \cdots \leq m_n.$ 

Since $m$ is rational in general, consider its floor (or integral part) $\lfloor m \rfloor$. 
By equation~\eqref{a-terms}, $m_2 = 2k_{22} \leq m,$ 
so there are at most $\lfloor\frac{m}{2}\rfloor$ values for 
$$
	m_2 = 2, 4, \ldots, \left\lfloor\frac{m}{2}\right\rfloor.
$$
Once again from equation~\eqref{a-terms} we obtain
\beq\label{m23}
	\frac{2k_{32}}{m_2} + \frac{2k_{33}}{m_3} = 1
\eeq
with $k_{33}\neq 0$,
implying
$$
	0 \leq k_{32}<\frac{m}{2}.
$$ 
If $\frac{m}{2} \not\in \Z,$ then there are $\lfloor\frac{m}{2}\rfloor+1$ choices for $k_{32},$ namely the integers from $0$ to $\lfloor\frac{m}{2}\rfloor.$ If $\frac{m}{2}\in \Z,$ then $0 \leq k_{32} \leq \frac{m}{2}-1,$ so there are $\frac{m}{2}$ choices for $k_{32}.$ 
In both cases, we get at most $\lfloor\frac{m}{2}\rfloor+1$ choices for $k_{32}.$ 
As for $k_{33},$ 
\eqref{m23} implies
$$0<2k_{33} \leq m_3 \le m,$$
so there are at most $\lfloor\frac{m}{2}\rfloor$ choices for $k_{33}.$ 
Since the choice of $k_{32}$ and $k_{33}$ determines $m_3,$ we have $\left(\lfloor\frac{m}{2}\rfloor+1\right)\lfloor\frac{m}{2}\rfloor$ choices for $m_3$, 
without accounting for different equations yielding the same solution. To determine the number of choice for $m_4,$ we use the equation $$\frac{2k_{42}}{m_2} + \frac{2k_{43}}{m_3}+ \frac{2k_{44}}{m_4} = 1.$$ The same analysis gives us $\lfloor\frac{m}{2}\rfloor+1$ choices for each of $k_{42}$ and $k_{43}$ and $\lfloor\frac{m}{2}\rfloor$ choices for $k_{44}$ due to the condition $k_{44}\neq 0.$ We thus have at most $\left(\lfloor\frac{m}{2}\rfloor+1\right)^2\lfloor\frac{m}{2}\rfloor$ choices for $m_4.$ In general, there are $\left(\lfloor\frac{m}{2}\rfloor+1\right)^{j-2}\lfloor\frac{m}{2}\rfloor$ choices for $m_j,$ where $2 \leq j \leq n.$ Altogether, we have obtained $$\left(\left\lfloor\frac{m}{2}\right\rfloor+1\right)^{0+1+ \cdots +(n-2)}\left(\left\lfloor\frac{m}{2}\right\rfloor\right)^{n-1}=\left(\left\lfloor\frac{m}{2}\right\rfloor+1\right)^{\frac{(n-2)(n-1)}{2}}\left(\left\lfloor\frac{m}{2}\right\rfloor\right)^{n-1}$$ possible multitypes  at $0$ of a pseudoconvex hypersurface $M\subset\C^n$ of finite $1$-type $m.$ This estimate significantly improves the one in \cite{ra}.

\ee

\section{Notation}
\label{notation}

For an $n$-tuple $z = (z_1, \ldots, z_n)$,
we shall use the short-hand notation 
  \[  
  	z_{[k,m]} := (z_k, z_{k+1}, \ldots, z_m),
	\quad
	1\le k \le m\le n. 
  \]
 We use the extended sets of nonnegative
 rationals and reals
 $$
      \1\Q_{\ge0} := \Q_{\ge0}\cup \{+\infty\},
      \quad
      \1\R_{\ge0} := \R_{\ge0}\cup \{+\infty\},
 $$
 and
  consider 
 real nonnegative  $n$-tuples of weights,
 or simply weights
 \beq\label{reg-wt}
 \mu= (1, \mu_2,  \ldots, \mu_n),
 \quad
 1=\mu_1> \mu_2\ge \ldots \ge \mu_n\ge 0,
 \quad
 \mu_i \in \R_{\ge0}.
 \eeq
 Following Catlin's notation \cite{catlinbdry},
 we also consider {\em inverse weights}
  $$
  	\Lambda= (1, \lambda_2, \lambda_3, \ldots, \lambda_n),
	\quad
	\lambda_i \in \1\Q_{\ge0}.
$$ 
For every weight $\mu= (1, \mu_2, \mu_3, \ldots, \mu_n),$ 
we have 
its associated inverse weight given by reciprocals
$$
	\Lambda= (1, \lambda_2, \ldots, \lambda_n)
	= (1, \mu_2^{-1}, \ldots, \mu_n^{-1})
$$ 
with the convention that $0^{-1} = +\infty$, $(+\infty)^{-1}=0$.
  Let 
  $$
  	(\alpha | \mu) := \alpha_1 \mu_1 + \alpha_2 \mu_2 + \cdots + \alpha_n \mu_n
$$ 
for $\alpha = (\alpha_1, \alpha_2, \dots, \alpha_n)$ a multi-index. 

Given a smooth real function $r(z,\bar z)$ defined 
in a neighborhood of $0$ in $\C^n$, 
a weight $\mu$ as in \eqref{reg-wt}
and nonnegative constant $C\ge0$
we write
$$
	r = O_\mu(C), 
	\quad \text{ resp. }
	r = o_\mu(C), 	  
$$
whenever 
$(\a+\b|\mu)\ge C$, 
resp.\
$(\a+\b|\mu)> C$
holds
for any 
nonzero monomial
$r_{\a\b}z^\a\bar z^\b$
in the Taylor expansion of $r$ at $0$.

Let 
$M\subset\C^n$
be an oriented smooth real hypersurface
defined in a neighborhood of a point $p=0$
by $r=0$ with $dr\ne 0$,
such that $r<0$ is the negative side
with respect to the orientation.
Recall that $M$
 is pseudoconvex if
and only if
the restriction of the complex Hessian of $r$
to the complex tangent space of $M$
is positive semidefinite.

We have the following elementary properties,
provided with short proofs for the reader's convenience.

\bl
	Let 
	$$
		r = -2 \re z_1 + f(z_{[2,n]}, \z_{[2,n]}),
	$$
	where $f$ is any smooth function.
	Then the domain given by $\{r<0\}$
	is pseudoconvex if and only if
	\beq\label{f}
		\sum_{j,k=2}^n f_{z_j\z_k} a_j \bar a_k \ge 0
	\eeq
	for all $(a_2, \ldots, a_n)\in \C^{n-1}$.
\el

\bpf
	The Levi form of the boundary $M := \{r=0\}$ is given
	by the restriction of
	the complex Hessian of $r$, given by the 
	left-hand side of \eqref{f},
	to the complex tangent subbundle of $M$.
	Since the latter projects surjectively to 
	$\{0\}\times \C^{n-1}$, the Levi form
	is positive semidefinite if and only if \eqref{f} holds,
	proving the statement.
\epf

\bl
	Given a weight $\mu$ as in \eqref{reg-wt},
	let 
	$$
		r = -2 \re z_1 + p(z_{[2,n]}, \z_{[2,n]}) + o_\mu(1),
	$$
	where $p$ is a weighted homogeneous
	polynomial of weight $1$.
	Assume that the domain given by $\{r<0\}$
	is pseudoconvex.
	Then the model domain given by 
	$\{r_0<0\}$, where
	$$
		r_0 := -2 \re z_1 + p(z_{[2,n]}, \z_{[2,n]}),
	$$
	is also pseudoconvex. \label{pscmodel}
\el

\bpf
	The proof is obtained by a simple weighted scaling argument.
	Consider the weighted dilation
	$$
		T_t(z_1, \ldots, z_n):= (z_1, t^{\mu_2} z_2, \ldots, t^{\mu_n} z_n).
	$$
	Then $r_0$ is invariant under composition with $T_t$, whereas
	for any function $f(z,\z)=o_\mu(1)$,
	the rescaled function $f(T_t(z), \1{T_t(z)})$
	 converges to $0$ uniformly on compacta as $t\to 0$.
	 The statement follows from the continuity 
	 of the complex tangent bundles and the
	  Levi form under the limit $t\to 0$.
\epf

Finally, we shall write
\beq\label{sim}
	A\sim B
\eeq
whenever there is a nonzero constant $c$
with $A=cB$.

\section{Catlin Multitype and Boundary Systems}
\label{bdrysystems}

We devote this section to defining the multitype notion
as Catlin introduced in \cite{catlinbdry} in order to characterize the vanishing order of the defining function in different directions.

In the previous section, we defined completely general weights and inverse weights, but Catlin restricts the inverse weights he considers to only those that could represent the vanishing of the defining function. 
We restrict the set of weights via two natural definitions:

\begin{defn}
An inverse weight $\Lambda= (1, \lambda_2, \lambda_3, \ldots, \lambda_n)$ is called {\it admissible} if for every $i,$ $1 \leq i \leq n,$ either $\lambda_i = + \infty$ or there exists a set of non-negative integers $a_1, \dots , a_i$ with $a_i >0$ such that $\displaystyle\sum_{j=1}^i {a_j}\lambda_j^{-1}=1.$ 
Let $\Gamma_n$ be the set of all admissible inverse weights ordered lexicographically. 
\end{defn}

\begin{defn}
Consider a smooth domain $\Omega \subset \C^n$ with defining function $r.$ 
An admissible inverse weight $\Lambda=(\lambda_1, \dots, \lambda_n) \in \Gamma_n$ is called {\it distinguished} at $z_0\in \d\Omega$
 if there exist holomorphic coordinates $(z_1, \dots, z_n)$ about $z_0$ with $z_0$ mapped to the origin such that if $\sum_{i=1}^n \: \frac{\alpha_i+ \bar \beta_i}{\lambda_i}<1,$ then $D^\alpha \bar D^{\bar \beta} r(0)=0,$ where $D^\alpha= \frac{\partial^{|\alpha|}}{\partial z^{\alpha_1}_1 \cdots \partial z^{\alpha_n}_n}$ and $\bar D^{\bar\beta}= \frac{\partial^{|\bar\beta|}}{\partial \bar z^{\bar\beta_1}_1 \cdots \partial \bar z^{\bar\beta_n}_n}.$ Let $\tilde \Gamma_n (z_0)$ be the set of distinguished weights at $z_0.$ \label{distinguishedweightdef}
\end{defn}

\begin{defn}
\label{multitypedef} The multitype $\mt(z_0)$ is defined to be lexicographically the smallest admissible weight $\mt(x_0) = (m_1, \dots, m_n)$ such that $\mt(z_0) \geq \Lambda$ for every admissible distinguished weight $\Lambda \in \tilde \Gamma_n (z_0).$
\end{defn}

\noindent Definitions~\ref{distinguishedweightdef} and~\ref{multitypedef} together prompt the following natural question:

\smallskip\noindent {\bf Question:} 
Let the multitype $\mt(z_0) = (m_1, \dots, m_n)$ be such that $m_n < + \infty.$ 
What are the multi-indices $\alpha$ and $\bar\beta$ satisfying $\sum_{i=1}^n \: \frac{\alpha_i+ \bar \beta_i}{\lambda_i}=1$ such that $D^\alpha \bar D^{\bar \beta} r(0)\neq0$
after a holomorphic change of variables mapping $z_0$ to the origin?

\smallskip\noindent Corollary~\ref{catlincase}, 
which we shall prove, 
not only gives an answer to this question but also identifies balanced terms in the defining function responsible for the condition $D^\alpha \bar D^{\bar \beta} r(0)\neq0.$ 

A priori, 
Definition~\ref{multitypedef} gives no indication how to compute the multitype $\mt(z_0).$ 
To achieve that, Catlin introduced in \cite{catlinbdry} the commutator multitype $\ct(z_0)$ computed by differentiating the Levi form
along certain lists of vector fields
arising from a geometric object called a {\em boundary system}.
He was then able to prove that this commutator multitype $\ct(z_0)$ equals the multitype
 for a pseudoconvex domain. 
 
 Recall from~\cite{catlinbdry} that
a boundary system is a collection of vector fields and real-valued functions $$\B_\nu (z_0)= \{r_1, r_{p+2}, \dots, r_\nu; L_2, \dots, L_\nu\}$$ for some $\nu \leq n$.
The first function in the boundary system is $r_1=r,$ the defining function. Let $p$ be the rank of the Levi form of $b \Omega$ at $z_0.$ Since $r=0$ defines a manifold, we can choose 
the vector field $L_1$ such that $L_1(r)=1.$ 
Recall from the beginning of this section that the multitype seeks to capture the vanishing order of the defining function in different directions. From the information we have so far, if $\ct(z_0)= (1, c_2, c_3, \ldots, c_n),$ then the first entry comes from the condition $L_1(r)=1,$ and subsequently, $c_2 = \dots = c_{p+1}=2,$ an entry of $2$ for every non-zero eigenvalue of the Levi form.  We choose vector fields of type $(1,0)$ $L_2, \dots, L_{p+1}$ such that $L_i(r)=\partial r (L_i)\equiv 0$ and the $p \times p$ Hermitian matrix $\partial \bar \partial r (L_i, L_j)(x_0)$ is nonsingular for $2 \leq i,j \leq p+1.$ We have kept Catlin's notation of round parentheses for the evaluation of forms on vector fields. 
If $p+1=n,$ our construction is finished; otherwise, we need to make sense of vanishing orders higher than two. Let us denote by $T^{\, (1,0)}_{p+2}$ the bundle composed of $(1,0)$ vector fields $L$ such that $\partial r (L)=0$ and $\partial \bar \partial r (L, \bar L_j)=0$ for $j=2, \dots, p+1.$ For $l \geq 3$ we denote by $\ls$ a list of vector fields $\ls= \{L^1, \dots, L^l\}$ and by $\ls  \partial r$ the function $$\ls  \partial r (z) = L^1 \cdots L^{l-2} \, \partial r \, ([L^{l-1},L^l]) (z)$$ for $z \in b \Omega.$ We are interested in lists $\ls$ such that $\ls  \partial r(z_0) \neq 0$ that are chosen in the most natural way possible. For every $j$ such that $p+2 \leq j \leq n,$ we will pick a $(1,0)$ vector field $L_j \in T^{\, (1,0)}_{p+2}$ and a corresponding real-valued function $r_j$ such that $L_j r_j \neq 0$ but $L_j r_i = 0$ whenever $i<j.$ The process is inductive. When $j = p+2,$ the simplest possible list $\ls$ of $(1,0)$ vectors in $T^{\, (1,0)}_{p+2}$ that can yield $\ls  \partial r(z_0) \neq 0$ consists of a smooth $(1,0)$ vector field $L_{p+2}$ and its conjugate $\bar L_{p+2}.$ If no such list exists, we set $c_{p+2}= \dots =c_n =+\infty,$ and we have finished our construction of the commutator multitype. If such a list exists, however, we choose a list $\ls_{p+2} = \{ L^1, \dots, L^l \}$ of minimal length $l$ and set $c_{p+2}=l.$ Let $\ls'_{p+2} = \{ L^2, \dots, L^l \}.$ 
Set $r_{p+2}= \re (\ls'_{p+2} \partial r)$ or $r_{p+2}= \im (\ls'_{p+2} \partial r)$ so that the condition $L_{p+2} r_{p+2} \neq 0$ holds. Define $S_{p+2}=\{L_{p+2}, \bar L_{p+2}\}.$

Now assume that for some integer $j-1$ with $p+2 \leq j-1<n$ we have already constructed finite positive numbers $c_1, \dots, c_{j-1}$ as well as real-valued functions $r_1,$ $r_{p+2}, \dots, r_{j-1};$ linearly independent smooth $(1,0)$ vector fields $L_2, \dots, L_{j-1};$ and lists $\ls_{p+2}, \dots, \ls_{j-1}$ such that the following properties hold:
\begin{enumerate}
\item $\ls_i \partial r(z_0) \neq 0$ for every $i,$ $p+2 \leq i \leq j-1;$
\item If $\ls_i = \{ L^1, \dots, L^l \},$ then $\ls'_i = \{ L^2, \dots, L^l \}$ and $r_i= \re (\ls'_i \partial r)$ or $r_i= \im (\ls'_i \partial r)$ in order that the condition $L_i r_i \neq 0$ holds;
\item $L_i r_k =0$ for $p+2 \leq k < i \leq j-1;$
\item Each of the lists $\ls_i = \{ L^1, \dots, L^l \}$ is
\begin{enumerate}[(a)]
\item i-{\it admissible} (in Catlin's terminology) meaning $L^1 \in S_i =\{L_i, \bar L_i\}$ and if $l_k$ is the number of times a vector from $S_k$ occurs in the list, then $$\displaystyle \sum_{k=p+2}^{i-1} \frac{l_k}{c_k} <1$$

and

\item {\it ordered} meaning $L^k \in S_{\alpha_k}$ for every $1 \leq k \leq l$ and $\alpha_1 \geq \alpha_2 \geq \dots \geq \alpha_l;$
\end{enumerate}
\item If $l^i_k$ equals the number of times $L_k$ and $\bar L_k$ occur in $\ls_i,$ then $l^i_k=0$ whenever $k>i$ and $\displaystyle \sum_{k=p+2}^{j-1} \frac{l^i_k}{c_k} = 1;$
\item All lists $\ls_i$ are of minimal length, namely if $\ls = \{ L^1, \dots, L^l \}$ is any ordered list, $l_k$ equals the number of times $L_k$ and $\bar L_k$ occur in $\ls,$ and $\displaystyle \sum_{k=p+2}^{j-1} \frac{l_k}{c_k} <1,$ then $\ls\partial r(z_0) = 0.$
\end{enumerate} 

We will show that we can chose a positive rational number $c_j$ such that properties (1)-(6) are fulfilled with $j$ replacing $j-1.$ Let $T^{\, (1,0)}_j$ denote the set of $(1,0)$ smooth vector fields $L$ such that $\partial \bar \partial r (L, \bar L_i)=0$ for $i=2, \dots, p+1$ and $L(r_k)=0$ for $k=1,p+2,\dots,j-1.$ 
For any smooth vector field $L_j \in T^{\, (1,0)}_j,$ we consider all ordered j-admissible lists $\ls.$ If for all such lists, $\ls\partial r(z_0) = 0,$ then we set $c_j = \dots = c_n = + \infty;$ otherwise, there exists at least one such list $\ls$ for which $\ls\partial r(z_0) \neq 0.$ We choose one of minimal length and denote it $\ls_j.$ The vector field $L_j$ used in its construction gets added to the collection $L_2, \dots, L_{j-1}.$ For $p+2\leq k \leq j$ let $l_k$ be the number of times a vector from $S_k=\{L_k, \bar L_k\}$ occurs in the list   $\ls_j,$ and let $c(\ls)$ denote the solution to the equation $$\sum_{k=p+2}^{j-1} \frac{l_k}{c_k}+\frac{l_j}{c(\ls)}=1.$$ Since $\ls_j$ is j-admissible, $c(\ls)\in \Q^+.$ Let $c_j = c(\ls)$ and let $r_j= Re (\ls'_j \partial r)$ or $r_j= Im (\ls'_j \partial r)$ so that the condition $L_j r_j \neq 0$ holds. All properties (1)-(6) are thus fulfilled. We continue this process until we have generated $\ct(z_0)= (1, c_2, c_3, \ldots, c_n),$ the {\it commutator multitype} at $z_0.$ Let $\nu \leq n$ be the highest index for which the entry $c_\nu$ is finite. The collection $$\B_\nu (x_0) = \{r_1, r_{p+2}, \dots, r_\nu; L_2, \dots, L_\nu\}$$ of functions and vector fields that we have generated in the process of computing $\ct(z_0)$ is called a {\it boundary system} of rank $p$ and codimension $n-\nu.$

\section{Estimates for Non-negative Homogeneous Polynomials}
\label{estimates}

\subsection{One Variable Case}
\begin{lem}
Let $P(z,\bar z)\ge0$ be a real
homogeneous polynomial of degree $2m$ in 
$\C$,
\[
P(z, \bar z) =
\sum_{k=-m}^m C_k z^{m+k} \bar z^{m-k}.
\]
Then 
$C_0 \ge 0$ and 
$ |C_{k}| \le  C_0$ for all $k$. 
Furthermore, if $P(z,\bar z) \not\equiv 0$, then $C_0 > 0.$
\label{estdegree1}
\end{lem}

\bpf
Considering the values of $z$ on $|z| = 1$, we observe that $\bar z = 1/z$, and
\[
P = \sum_{k=-m}^m C_k z^{2k} \ge 0, \quad |z| = 1.
\]
We parametrize $|z|=1$ as $\gamma (\theta) = e^{i \theta}$ on $[-\pi, \pi]$ and observe that
\[
0 \le \frac{1}{2\pi } \int_{-\pi}^\pi P(\theta) d \theta = C_0,
\]
whence $C_0 \ge 0$. Note that if $P > 0$ for some $z$ in a neighbourhood of $0,$ then $\displaystyle \int_{-\pi}^\pi P(\theta) d \theta >0$ so $C_0 > 0.$

Since for any $k \neq 0$ and $z$ with $|z|=1$,
\[
0 \le (2\re z^k)^2 = (z^k + \bar z^k)^2 = (z^k + 1/z^k)^2 = z^{2k} + z^{-2k} + 2
\]
and
\[
0 \ge (2i\im z^k)^2 = (z^k - \bar z^k)^2 = (z^k - 1/z^k)^2 = z^{2k} + z^{-2k} - 2,
\]
we have
\[
0 \le \frac{1}{2\pi } \int_{-\pi}^\pi (e^{i\theta k} + e^{-i\theta k})^2 P(\theta) d \theta = C_{-k} + C_{k} + 2C_{0}
\]
and
\[
0 \ge \frac{1}{2\pi }\int_{-\pi}^\pi (e^{i\theta k} - e^{-i\theta k})^2 P(\theta) d \theta = C_{-k} + C_{k} - 2C_{0},
\]
which immediately gives $ |C_{-k} + C_{k}| \le 2C_{0} $. Since $C_{k} = \bar C_{-k}$, we have
\begin{equation}
\label{eq:real_part_is_smaller}
|\re C_{-k} | = |\re C_{k} | \le C_{0}.
\end{equation}

Furthermore, by rotating $z$ we may assume $C_k\in\R$,
and hence the above inequality yields $|C_k|\le C_0$ as desired.
\epf

\medskip
Inspired by these estimates, we seek to divide the terms of a non-negative homogenous polynomial $P(z, \bar z)$ in $\C^n$ into the terms that control others and those that are controlled. 

\begin{defn}
Let a universal constant $M>0$ be given, and let $P(z,\bar z)$ be a non-negative homogeneous polynomial of degree $2m$ in $\C^n.$ If 
\[ 
	P(z, \bar z) = \sum_{|\alpha|+|\bar\beta|=2m} C_{\alpha \bar\beta} \, z^\alpha {\bar z}^{\bar \beta}, 
\] 
then a coefficient $C_{\alpha \bar\beta}$ is called $M$-dominant if $|C_{\alpha' \bar\beta'}| \le M\, |C_{\alpha \bar\beta}|$ 
for all $\alpha', \bar\beta'$ such that $|\alpha'|+|\bar\beta'|=2m.$
\end{defn}

\subsection{Newton Polygon Lemma}

\begin{lem}\label{newton}
Let $P(x,y)$ be a non-negative homogeneous polynomial of degree $2m$ in $\R^p \times \R^q$ for $x \in \R^p$ and $y \in \R^q.$ If \[ P(x, y) = \sum_{p+q=2m} P_{pq}, \] where $P_{pq}(\lambda x, \mu y)= \lambda^p \mu^q P(x,y),$ then $P_{p_0 q_0}\geq 0$ if either $p_0 = \max p$ or $q_0 = \max q.$ \label{newtonpolygon}
\end{lem}

\smallskip\noindent {\bf Proof:} We use a scaling argument. $P(x,y) \geq 0$ for all $(x,y) \in \R^p \times \R^q$ implies for any $t \in \R,$ $t \neq 0,$ we have that $P(t x, t^{-1} y) \geq 0.$ Letting $t \to \infty$ shows $P_{p_0 q_0}\geq 0$ when $p_0 = \max p.$ Letting $t \to 0$ shows $P_{p_0 q_0}\geq 0$ when $q_0 = \max q.$ \qed

\smallskip

\begin{rem}
The statement of
Lemma~\ref{newtonpolygon} also holds for weighted homogenous polynomials.
\end{rem}

\section{Normal forms in the pseudoconvex case}
\label{normalforms}

\subsection{Setup}
Consider a smooth real hypersurface $M$ passing through $0$,
defined near $0$ by $r=0$, where
\[
r = -2\re z_1 + p(z_{[2,n]},\bar z_{[2,n]}) + 
o_\mu(1),
\]
$p$ is a weighted homogeneous polynomial 
of weight $1$,
and
where the components $z_1,\ldots, z_n$ and their conjugates
are assigned the corresponding weights 
$1, \mu_2, \mu_3, \ldots, \mu_n$.
In other words, we can write
$$
	\displaystyle p(z_{[2,n]},\bar z_{[2,n]}) = 
	\sum_{(\alpha +  \beta | \mu) =1} 
	C_{\alpha  \beta} z^\alpha \bar z^{ \beta}.
$$
Together with $M,$ we consider its model hypersurface $M_0$ defined by the weighted homogeneous part
\[
r_0=0,
\quad r_0: = -2\re z_1 + p(z_{[2,n]}, \bar z_{[2,n]}).
\]
Note that by elementary scaling argument, if $M$ is pseudoconvex, then so is its model $ M_0$ as seen in Lemma~\ref{pscmodel}.

As customary, we shall assume that $p(z_{[2,n]}, \bar z_{[2,n]})$ 
does not contain any pure (harmonic) monomials
 of the form $z^\alpha$ or $\bar z^\beta$.
 Otherwise, they can be always eliminated by
 a holomorphic transformation
 \[
 	(z_1, z_{[2,n]}) \mapsto (z_1 + h(z_{[2,n]}),  z_{[2,n]}).
\]

In what follows, homogeneity will be gauged with respect to a regular weight $\mu$

\subsection{First Step}
\begin{lem}\label{5.1}
Let $M_0\subset \C^n$, $n\ge 2$, be a pseudoconvex 
model hypersurface
defined by
\[
	r_0 =0, \quad  r_0= -2\re z_1 + p(z_{[2,n]}, \bar z_{[2,n]}),
\]
where $p$ is a weighted homogeneous polynomial of weight $1$ 
with respect to some weights
\[
	\mu= (\mu_1, \mu_2, \mu_3, \ldots, \mu_n), 
	\quad
	1=\mu_1>\mu_2\ge \mu_3 \ge \ldots \ge \mu_n \ge0.
\]
Let $s$ be such that
$$
	\mu_2 = \ldots = \mu_s  > \mu_{s+1},
$$
and assume that
\beq\label{p2r}
	p_{[2,s]}(z_{[2,s]}, \bar z_{[2,s]}) :=
	p(z_{[2,s]},0, \bar z_{[2,s]}, 0) \not\equiv 0.
\eeq

Then after a unitary change of the variables $(z_2, \dots, z_s)$, 
we have
\beq\label{p2-def}
	p_2(z_2, \bar z_2) := p(z_2,0, \bar z_2, 0) \not\equiv 0,
\eeq
which is a plurisubharmonic homogeneous polynomial
of even degree $2 k_{22}$, where $k_{22} := \frac{1}{2 \mu_2}$.
Furthermore, $p_2$ has the form
\beq \label{p2}
	p_2 = \sum_{j=-k_{22}+1}^{k_{22}-1} C_{2,j} \, z_2^{k_{22}+j} \bar z_2^{k_{22}-j}, 
\eeq
where $C_{2,0}>0
$ is $\frac{1}{2 \mu_2}$-dominant among all the coefficients of $p_2$.
\label{firststeprotation}
\end{lem}

\bpf
By our assumptions, $p_{[2,s]}$ is a nonzero
weighted homogenous polynomial
in $(z_{[2,s]}, \bar z_{[2,s]})$.
Then clearly after a unitary transformation of $z_{[2,s]}$,
we may assume that \eqref{p2-def} holds with $p_2$
of the form \eqref{p2}.
%

Using pseudoconvexity of $M_0$ 
 we obtain
$$ 
	\partial_{z_2 \bar z_2} p_2 \ge 0,
$$ 
where $p_2$ is defined by \eqref{p2-def}.
Since $p_2$ does not contain any harmonic terms by our assumption, 
each monomial in $(z_2, \z_2)$
will contribute to $\partial_{z_2 \bar z_2} p_2$.
We compute
\[ 
	\partial_{z_2 \bar z_2} p_2 = \sum_{j=-k_{22}+1}^{k_{22}-1} C_{2,j}(k_{22}+j)(k_{22}-j) \, z_2^{k_{22}+j-1} \bar z_2^{k_{22}-j-1} .
\]
By Lemma~\ref{estdegree1} applied to $\partial_{z_2 \bar z_2} p_2,$ 
we conclude 
$C_{2,0} > 0,$
and

\[
	|C_{2,j}| \leq 
	\frac{k_{22}^2}
	{k_{22}^2-j^2} C_{2,0}
	\le
	\frac{k_{22}^2}
	{2k_{22}-1} C_{2,0}
	<
	k_{22} C_{2,0}
	 ,
	 \quad
	 j \ne 0. 
\]

 Note that $C_{2,0}$ is the coefficient of the balanced term $C_{2,0} \, z_2^{k_{22}} \bar z_2^{k_{22}}$ in $p_2.$ Since $r_0$ has weight $1$ with respect to the weight $\mu,$ we have
 $
 	k_{22} =\frac{1}{ 2\mu_2}
$
and $2 k_{22}$ is the degree of $p_2$. 
\epf

\bigskip
We shall now apply Lemma~\ref{5.1}
with the weights given by the Catlin multitype.
Using the previous lemma, we can normalize the first non-trivial function 
in the Catlin boundary system:

\begin{lem}
Let $M$ be a pseudoconvex hypersurface 
defined by $r=0$, where
\[
	r =
	 -2\re z_1 + p(z_{[2,n]}, \bar z_{[2,n]})
	 + o_\mu(1),
\]
such that $p$ is a weighted homogeneous polynomial of weight $1$
with respect to the Catlin multitype at $0$,
\[
\Lambda= (1, \lambda_2, \lambda_3, \ldots, \lambda_n).
\]
Assume the Levi rank at $0$ is $s_0$ and $2<\lambda_{s_0+2} = \dots = \lambda_{s_0+s_1+1}< \lambda_{s_0+s_1+2}< +\infty.$ Let $$\B_n (0) = \{r_1, r_{s_0+2},\dots, r_n; L_2, \dots, L_n\}$$ be any boundary system at $0.$ Then there
exists a holomorphic change of coordinates at $0$ 
 preserving the multitype so that the boundary system in the new coordinates $$\tilde \B_n (0) = \{\tilde r_1,\tilde  r_{s_0+2}, \dots, \tilde r_n; \tilde L_2, \dots, \tilde L_n\},$$ has the simplest possible first functions $r_{s_0+2}, \dots, r_{s_0+s_1+1}$ given by
\beq
	r_j = \re z_j+o\left(\frac{1}{\lambda_j}\right)
\eeq
\label{firststepfixr2}
for $s_0+2 \leq j \leq s_0+s_1+1.$
\end{lem}


\bpf By a Chern-Moser type argument \cite{CM}, we may assume that at $0$ the Levi rank $s_0=0.$ 
Furthermore, it is easy to see this lemma reduces to proving that at the level of the model hypersurface $r_0=0,$ where \[
	r_0 =
	 -2\re z_1 + p(z_{[2,n]}, \bar z_{[2,n]}),
\] for $s_0+2 \leq j \leq s_0+s_1+1$ we can bring $r_j$ to the form $r_j=\re z_j$ via a holomorphic polynomial change of the variables.
We start by noting that 
the assumption \eqref{p2r}
must hold for the Catlin's multitype of the model hypersurface, which is the same as that of the original hypersurface.
Indeed, otherwise we could increase the inverse weights
$$
	\lambda_2 =\ldots = \lambda_{s_1+1}
$$
and decrease the remaining bigger inverse weights,
still keeping the total weighted degrees of
the terms in $p$ greater or equal $1$, which would contradict that assumption that $(1, \lambda_2, \lambda_3, \ldots, \lambda_n)$ is the Catlin multitype at $0$. For the moment, assume $s_1=1.$

We can apply Lemma~\ref{firststeprotation} to obtain
a decomposition
\[ 
	r_0 = -2\re z_1 + p_2(z_2, \bar z_2) + q_2(z_{[2,n]}, \bar z_{[2,n]}),
\]
where $p_2$ satisfies \eqref{p2} and
$$
	q_2(z_2,0,\bar z_2,0) \equiv 0.
$$
We would like to make another change of variables that would ensure the function $r_2$ in the boundary system has the simplest possible expression, namely $r_2 = \re z_2.$ Note that the boundary system contains a function $r_2$ due to the assumption $2<\lambda_2 < +\infty.$ To determine what change of variables needs to be made, we notice that  regardless of what form $r_2$ initially assumes, after the change of variables mandated by the application of Lemma~\ref{firststeprotation},
\begin{equation}
\label{rp+2prep}
\frac{1}{(k'-1)! k''!}
\partial^{k' - 1}_{z_2} \partial^{k''}_{\bar z_2} r_0
= C_{2,0} \, z_2 + C_{2,-1} \, \bar z_2
+ T(z_{[3,n]},\bar z_{[3,n]} ),
\end{equation}
where $k'+k''=2 k_{22}=\lambda_{2}$ and $ T(z_{[3,n]},\bar z_{[3,n]} )$ is the sum of terms coming from 
differentiation of $q_2$ in the expression for $r_0.$ Note that $T$ cannot depend on $z_2$ due to the weight restriction.

It is beneficial to split $q_2$ as 
\begin{equation*}
q_2 =
\sum_{k+l = 2k_{22}-1} z_2^k \bar z_2^l T_{kl}(z_{[3,n]}, \bar z_{[3,n]})\\
%
+
S_2(z_{[2,n]}, \bar z_{[2,n]}),
\end{equation*}
where all monomials in 
$S_2$ have their total degree in $(z_2, \bar z_2)$
less than $2k_{22}-1$.

By the pseudoconvexity of $M_0$, we see that for $j \ge 3$,
\[
\partial_{z_j \bar z_j} r_0 =
	\sum_{k+l = 2k_{22}-1} z_2^k \bar z_2^l  
	\partial_{z_j \bar z_j}T_{kl}(z_{[3,n]}, \bar z_{[3,n]})
+ \partial_{z_j \bar z_j} S_2 \ge 0.
\]
Choosing
$|z_2| >> |z_j|$ for all $j \ge 3$,
 we conclude
\[
	\sum_{k+l = 2k_{22}-1} z_2^k \bar z_2^l  
	\partial_{z_j \bar z_j}T_{kl}(z_{[3,n]}, \bar z_{[3,n]})
	\ge 0, 
	\quad
	j \ge 3.
\]
Since the left-hand side changes the sign when $z_2$ does, we obtain
%
\[
\partial_{z_j \bar z_j} T_{kl} = 0, \ j \ge 3.
\]
Since the $z_j$-direction can be rotated 
arbitrarily by a change of coordinates,
by a similar argument, the whole Levi form of each $T_{kl}$
must vanish identically,
which means that each $T_{kl}$ is harmonic,
i.e.\ a sum of holomorphic and antiholomorphic functions 
in $z_{[3,n]}$.

In particular, after the change of variables mandated by the application of Lemma~\ref{firststeprotation},
in the boundary system
$$
	L_2:= \d_{z_2} + p_{z_2} \d_{z_1},
	\quad
	\ls = \{\bar L_2, L_2, \bar L_2, \ldots, L_2, \bar L_2\},
$$
where $L_2$ and $\bar L_2$ appear $k'-1$ and $k''$ times respectively. Note that if $\lambda_2 = \lambda_3 = \dots = \lambda_{s_1+1},$ another linear change of variables to gather all terms coming from $z_3, \dots, z_{s_1+1}$ into $z_2$ may be required in order to achieve the expression for $L_2$ claimed above.

Since $p$ does not depend on $z_1$,
we can ignore differentiations in that direction
in the expressions of $\ls \d r$.
%
Hence we can compute $r_2$ up to a constant as
$$
	r_2 = \ls \d r  \sim   \re \partial^{k' - 1}_{z_2} \partial^{k''}_{\bar z_2} r_0
	 \sim \re( C_{2,0} \, z_2 + C_{2,-1} \, \bar z_2 + T),
$$
where
\[
T = \phi_2(z_{[3,n]}) + \overline{\psi_2(z_{[3,n]})}
\]
with $\phi_2$ and $\psi_2$ holomorphic. 
By a rotation in $z_2$, we can assume 
$C_{2,-1}$ real and nonnegative.
Since $C_{2,0}>0$,
we can consider the change of variables,
\begin{equation}
\label{eq:change_of_vars_z2}
z_2 = \frac{1}{(C_{2,0} + \bar C_{2,-1})} z'_2 - (\phi_2 + \psi_2)
\end{equation}
leading to
$r_2 \sim \re z_2$ and finally by scaling to
\[ r_2 = \re z_2.  \]
Now let $s_1>1.$ For every $j$ such that $3 \leq j \leq s_1+1,$ we carry out the same procedure noting that due to weight considerations the change of variables required to transform $r_j$  into $r_j = \re z_j$ will not affect variables $z_1, \dots, z_{j-1}.$
\epf

\smallskip
\begin{rem}
A similar argument cannot be carried out to normalize
the next boundary system function $r_{s_0+s_1+2}$ if $\lambda_{s_0+s_1+1} < \lambda_{s_0+s_1+2}.$ A counterexample is given at the end of the paper in Section~\ref{r3counterexample}. 
\end{rem}

\medskip\noindent {\bf Proof of Theorem~\ref{flatnessprop}:} At the level of the model hypersurface, $r_0=0,$ with \[
	r_0 =
	 -2\re z_1 + p(z_{[2,n]}, \bar z_{[2,n]}),
\] Lemma~\ref{firststepfixr2} shows that after a change of variables $\tilde r_j = \re z_j$ for $s_0+2 \leq j \leq s_0+ s_1+1.$ Therefore, $\tilde L_j=\frac{\partial}{\partial z_j}$ and $
\tilde L_{s_0+ s_1+2} =\frac{\partial}{\partial z_{s_0+ s_1+1}}+o\left(-\frac{1}{\lambda_{s_0+ s_1+2}}\right)$ for the model hypersurface as $\tilde L_{s_0+ s_1+2}$ is chosen so that $\tilde L_{s_0+ s_1+2} \tilde r_{s_0+ s_1+1} = 0$, and by Catlin own normalization result in \cite{catlinbdry}, Proposition 5.3, $$r_0=2\re z_1+\sum_{j=2}^{s_0+1} |z_j|^2+f_1(z_{s_0+2}, \dots, z_n)+ 2 \re \left(\sum_{j=2}^{s_0+1} z_j f_j(z_{s_0+2}, \dots, z_n)\right).$$ The proposition follows.
\qed

\medskip\noindent {\bf Proof of Corollary~\ref{C3case}:} We apply Proposition~\ref{flatnessprop} with $s_0=0$ and $n=3.$
\qed

\subsection{Second Step}

\begin{lem}
Let $M_0$ be a pseudoconvex hypersurface with polynomial defining function 
\[
r_0 = -2\re z_1 + p(z_{[2,n]}, \bar z_{[2,n]}).
\]
such that $p$ is a weighted homogeneous polynomial in $z_2, \dots, z_n$
 of total weight $1$ with respect to the weights
\[
\mu= (1, \mu_2, \mu_3, \ldots, \mu_n).
\]
Let $s$ be such that
$$
	\mu_3 = \ldots = \mu_s > \mu_{s+1},
$$
and assume that 
\begin{equation}
\label{eq:nondegeneracy_z3}
	p(z_{[2,s]},0, \bar z_{[2,s]},0) \not\equiv p_2(z_2, \bar z_2),
\end{equation}
where $p_2$ is given by \eqref{p2-def}.

Then after a unitary change of variables $z_3, \dots, z_s$,
the polynomial
$p$ admits a decomposition
\begin{equation}
\label{eq:decomposition_z3}
	p=p_2(z_2, \bar z_2) + p_3(z_{[2,3]}, \bar z_{[2,3]}) + q_3(z_{[2,n]}, \bar z_{[2,n]}),
\end{equation}
where $p_3$ is of degrees $2k_{23}$ and $2k_{33}$
in $(z_2,\bar z_2)$ and $(z_3, \bar z_3)$
respectively, $q_3$ has only terms of degree less than $2k_{33}$
in $(z_3, \bar z_3)$,
and $p_3$ contains a non-zero term 
$$
	C |z_2|^{2k_{23}} |z_3|^{2k_{33}} \ge 0, 
	\quad k_{33}>0.
$$

\label{secondsteprotation}
\end{lem}

\smallskip\noindent {\bf Proof:} 
Just like at the beginning of the proof of Lemma~\ref{firststepfixr2}, we write $r_0$ as
\[
r_0 = -2\re z_1 + p_2(z_2, \bar z_2) + q_2(z_{[2,n]}, \bar z_{[2,n]} ).
\]
After a possible unitary change of variables in $z_{[3,s]}$, 
we can assume that 
$$q_2(z_{[2,3]},0, \bar z_{[2,3]},0)\not\equiv 0.$$ 
Using Lemma~\ref{newtonpolygon} for guidance, we next identify the non-zero terms in $q_2$ 
of the highest (total) degree $d_3$ in $(z_3, \bar z_3)$ 
and denote their sum by 
$p_3(z_{[2,3]}, \bar z_{[2,3]})$.
Then $p$ can be decomposed as
$$
	p =  p_2(z_2, \bar z_2)
	+
	p_3(z_{[2,3]}, \bar z_{[2,3]})
	+ q_3(z_{[2,n]}, \bar z_{[2,n]} ),
$$
where 
all monomials in 
$q_3(z_{[2,3]}, 0, \bar z_{[2,3]}, 0)$
have degree less than $d_3$ in $(z_3, \bar z_3)$.

We shall consider the inequality
\beq\label{diag}
	(\d_{z_2} + t\d_{z_3}) (\d_{\bar z_2} + t\d_{\bar z_3}) 
	(p_2(z_2, \bar z_2) + p_3(z_{[2,3]}, \bar z_{[2,3]}) + q_3(z_{[2,3]},0, \bar z_{[2,3]},0)) \ge 0,
\eeq
that follows from the pseudoconvexity of $M_0$,
where $t$ is an arbitrary parameter.
Identifying terms of the highest degree in $(z_3, \bar z_3)$ we obtain
$$
	\d_{z_2\bar z_2} p_3 \ge 0.
$$

We first assume that 
\beq\label{z22}
	\d_{z_2\bar z_2} p_3 \equiv 0,
\eeq
i.e.\ all terms in $p_3$ are harmonic in $z_2$.
Then again identifying terms of the highest degree in $(z_3, \bar z_3)$ 
under this assumption, we obtain
$$
	2t \re \d_{z_2\bar z_3} p_3 + \d_{z_2\bar z_2} \tilde q_3 \ge 0,
$$
where $\tilde q_3$ is the sum of certain terms from $q_3$.
Since $t$ is arbitrary, we must have
$$
	\d_{z_2\bar z_3} p_3 \equiv 0.
$$
Since $p_3$ has no harmonic terms and all terms are harmonic in $z_2$,
the only possibility remaining is that $p_3$ is independent of $z_2$.
But then, since $p_3$ is nonzero and has no harmonic terms, we must have
\beq\label{d33}
	\d_{z_3\bar z_3} p_3 \not\equiv 0.
\eeq

On the other hand, if \eqref{z22} does not hold,
we obtain the polynomial $\d_{z_2\bar z_2} p_3\ge 0$,
which for any generic fixed $z_2$, is non-constant and homogeneous in $(z_3,\bar z_3)$,
which again implies \eqref{d33}.

Thus in all cases, we must have \eqref{d33}.
Applying Lemma~\ref{estdegree1}
to $\d_{z_3\bar z_3} p_3$ for fixed $z_2$ (when it does not identically vanish),
we conclude that the degree $d_3$ in $(z_3, \bar z_3)$
is even, $d_3=2k_{33}$, and $\d_{z_3\bar z_3} p_3$ contains 
nonzero terms of the form 
$$\tilde p(z_2,\bar z_2)z_3^{k_{33}}\bar z_3^{k_{33}} \ge 0.$$
Applying again Lemma~\ref{estdegree1}, this time to $\tilde p$,
we conclude that $p_3$ contains a nonzero term
$$
	C z_2^{k_{23}} z_3^{k_{33}} \bar z_2^{k_{23}}  \bar z_3^{k_{33}} \ge 0,
$$
as desired.
\qed

\subsection{Inductive Step}
The general inductive step 
will be obtained from the following result.

\begin{lem}
Let $M_0$ be a pseudoconvex (model) hypersurface through $0$ with polynomial defining function 
\[
r = -2\re z_1 + p(z_{[2,n]}, \bar z_{[2,n]}).
\]
such that $p$ is a weighted homogeneous polynomial in $z_2, \dots, z_n$ of weight $1$ 
with respect  to the weights
\[
\mu= (1, \mu_2, \ldots, \mu_n),
\quad
1>\mu_2\ge\ldots\ge \mu_n.
\]
Assume that we have already shown that
\[
\begin{aligned}
p(z_{[2,n]}, \bar z_{[2,n]}) &= 
p_2(z_2, \bar z_2)
+
p_3(z_{[2,3]}, \bar z_{[2,3]})
+
\cdots
\\
&
+
p_{m-1}(z_{[2,m-1]},\bar z_{[2,m-1]})
+ q_{m-1}(z_{[2,n]}, \bar z_{[2,n]}).
\end{aligned}
\]
where
$$
	q_{m-1}(z_{[2,m-1]},0,\bar z_{[2,m-1]},0) \equiv 0.
$$
Let $s$ be such that
$$
	\mu_m = \ldots = \mu_s  > \mu_{s+1},
$$
and assume that 
\begin{equation}
\label{eq:nondegeneracy_z_n}
	q_{m-1}(0,z_{[m,s]},0, 0,\bar z_{[m,s]},0) \not\equiv 0.
\end{equation}

Then after a unitary change of the variables $(z_m, \dots, z_s)$, 
$q_{m-1}$ admits the decomposition
\begin{equation}
\label{eq:decomposition_z_n}
	q_{m-1} = p_m(z_{[2,m]}, \bar z_{[2,m]}) + q_m(z_{[2,n]}, \bar z_{[2,n]}),
\end{equation}
where 
$p_m$ is a homogeneous polynomial of degree $2k_{mm}>0$ in $(z_m, \bar z_m)$
whose expansion contains a term 
%
\[
C |z_2|^{2 k_{2m}} |z_3|^{2 k_{3m}} \cdots |z_m|^{2 k_{mm}},
\quad C > 0,
\]
and
$q_m$ has only terms of degree less than $2k_{mm}$
in $(z_m, \bar z_m)$.
\label{inductivesteprotation}
\end{lem}

\smallskip\noindent {\bf Proof:}
After a possible unitary change of variables in $z_{[m,s]}$, 
we can assume that 
$$q_{m-1}(0,z_m,0, 0,\bar z_m,0)\not\equiv 0.$$ 
We next identify the non-zero terms in $q_{m-1}$ 
of the highest degree $d_m>0$ in $(z_m, \bar z_m)$ 
and denote their sum by 
$p_m(z_{[2,m]}, \bar z_{[2,m]})$.
We thus write
\[
	q_{m-1} = p_m(z_{[2,m]}, \bar z_{[2,m]}) + q_m(z_{[2,n]}, \bar z_{[2,n]}),
%
\] 
where all monomials in $q_{m}(z_{[2,m]},0, \bar z_{[2,m]},0)$ have degree less than $d_m$ in $(z_m, \bar z_m)$ by construction.

We shall now use the pseudoconvexity assumption on $M_0.$ 
For any $j,$ $2\leq j<m$, 
and any arbitrary real parameter $t,$ consider
\beq\label{diagm}
\begin{aligned}
	(\d_{z_j} + t\d_{z_m}) (\d_{\bar z_j} + t\d_{\bar z_m}) &
	\Big(
		p_2(z_2, \bar z_2) + p_3(z_{[2,3]}, \bar z_{[2,3]}) + \cdots \\&+
p_{m-1}(z_{[2,m-1]},\bar z_{[2,m-1]})
+
p_{m}(z_{[2,m]},\bar z_{[2,m]})\\&
+
q_{m}(z_{[2,n]}, \bar z_{[2,n]})
\Big) \ge 0.
\end{aligned}
\eeq
Identifying the highest degree terms in $(z_m, \bar z_m),$ we obtain that
\beq\label{pos}
\partial_{z_j\bar z_j} p_m  \ge 0,
\quad j<m.
\eeq

We first assume $\partial_{z_j\bar z_j} p_m  \equiv 0$ for all $j<m$,
i.e.\ all terms of $p_m$ are harmonic in $z_j$.
 Then looking at the highest degree terms in $(z_m, \bar z_m)$ in \eqref{diagm} yields 
$$
	2t \re \d_{z_j\bar z_m} p_m + \d_{z_j\bar z_j} \tilde q_m \ge 0,
$$
where $\tilde q_m$ consists of the sum of the terms of degree $d_m-1$ 
in $(z_m, \bar z_m)$
from $p_2+\cdots+ p_m + q_m$.
Given that $t$ is arbitrary, we conclude
$$
	\d_{z_j\bar z_m} p_m \equiv 0.
$$
Note that $p_m$ contains no harmonic terms and by our assumption, all terms of $p_m$ are harmonic in $z_j$ for all $j<m$.
Hence any nonzero term of $p_m$ must have both $z_m$ and $\bar z_m$,
that is
we must have that 
\beq\label{mmDpm}
\partial_{z_m\bar z_m} p_m  \not\equiv 0.
\eeq

Now, on the contrary, assume that $\partial_{z_j\bar z_j} p_m  \not\equiv 0$
for some $j<m$.
For any generic fixed $z_2, \dots, z_{m-1}$,
 the polynomial $\partial_{z_j\bar z_j} p_m$ is non-constant, non-negative by \eqref{pos}, 
 and homogeneous in $(z_m, \bar z_m)$ of degree $d_m>0$. 
 Clearly, \eqref{mmDpm} must hold in this case as well. Therefore, regardless of the case, \eqref{mmDpm} holds. 
 
 We claim that by Lemma~\ref{estdegree1} inductively applied to $\partial_{z_m\bar z_m} p_m,$ the expansion of the polynomial $\partial_{z_m\bar z_m} p_m$ contains a term
\beq\label{c-term}
	C |z_2|^{2 \alpha_2}|z_3|^{2 \alpha_3} \cdots |z_m|^{2 \alpha_m},
	\quad \alpha_m>0, 
	\quad C> 0.
\eeq
Indeed, first keep $z_2, \dots, z_{m-1}$ fixed and apply Lemma~\ref{estdegree1} to $\partial_{z_m\bar z_m} p_m,$ a non-constant, non-negative, 
 and homogeneous polynomial in $(z_m, \bar z_m)$ of degree $d_m-2$. 
We conclude that the sum of the terms in $p_m$ having equal degrees in $z_m$ and $\bar z_m$ is of the form
  $$
 	P^{m-1}(z_2, \dots, z_{m-1}, \bar z_2, \dots, \bar z_{m-1})  |z_m|^{2 \alpha_m} \geq 0,
	\quad
	\alpha_m = d_m/2 > 0,
 $$
 where 
 $P^{m-1}$ is a nonzero weighted homogeneous polynomial.

Let $l$ be the highest index among $2, \dots, m-1$ for which $P^{m-1}$ has  degree $d_l >0$ in $(z_l, \bar z_l).$ If no such $l$ exists, then $P^{m-1}$ is constant and positive, and we are done; otherwise we extract the sum $\tilde P^{m-1}$ of terms of the top degree $d_l$ in $(z_l, \bar z_l)$, which is not identically zero, and nonnegative in view of Lemma~\ref{newton},
 keep $z_2, \dots, z_{l-1}$ fixed and apply Lemma~\ref{estdegree1} to $\tilde P^{m-1}$ viewed as a homogeneous polynomial in $(z_l, \bar z_l).$ 
 Proceeding inductively, we see that 
 $
 p_m$ contains a non-zero term \eqref{c-term}
maximizing the multidegree in $(z_2,\z_2), \ldots, (z_m, \z_m)$
in the reversed lexicographic order 
as claimed.

\qed

\medskip\noindent {\bf Proof of Theorem~\ref{main}:} 
Our Main Theorem is a consequence of Lemmas~\ref{firststeprotation}, \ref{secondsteprotation}, and \ref{inductivesteprotation}.  
In fact, at each step, either the nonvanishing assumption in the Lemmas holds, and hence 
we obtain a positive $A_j$ or we can lower the weight $\mu$ lexicographically starting with 
$\mu_j$. In the latter case, either we regain the nonvanishing assumption for a lower $\mu_j$
and keep applying the Lemmas with lower weights, or no further term is left and the proof is complete.
\qed

\medskip\noindent {\bf Proof of Corollary~\ref{catlincase}:} This result follows from Theorem~\ref{main}. \qed

\section{Counterexample to a boundary system normalization}
\label{r3counterexample}

We would like to show via an example that the kind of normalization of function $r_2$ in Catlin's boundary system that we carried out in Theorem~\ref{flatnessprop} fails for $r_3.$ Let the defining function of the domain be given by 
\[
	r_0 = -2\re z_1 + p(z_{[2,4]}, \bar z_{[2,4]}),
\]
where $p$ is a weighted homogeneous polynomial chosen so that the weight of $z_2^2$ equals the weight of $z_3^3.$ For example, the polynomial
\[
\begin{split}	
p(z_{[2,4]}, \bar z_{[2,4]}) &= |z_2|^6 
	+ |z_2|^2 |z_3|^6 
	+ |z_2|^4 |z_3|^2 |z_4|^2
	+|z_2|^2|z_3|^4 |z_4|^4\\&
	+ 2\epsilon \re ( |z_2|^2 z_3^2 \bar z_3^3  |z_4|^2)
	+ |z_3|^8|z_4|^2
\end{split}
\]

is homogeneous of weight $1$ with respect to

$$
	\Lambda = \left(1, \frac16, \frac19, \frac1{18}\right).
$$
Note that
$$
	f:=\d_{z_2}\d_{\bar z_2} \d_{z_3}^2 \d_{\bar z_3}^3 p
	= c_1 z_3 + c_2 |z_4|^2,
$$
and hence $r_3=\re f$
cannot be tranformed into $c\re z_3$
by any holomorphic coordinate change.

We would like to show that $p$ is plurisubharmonic when $\epsilon$ is small, where $0<\epsilon<1.$
First, we observe that if $z_2=0,$ $z_3=0,$ or $z_4=0,$ the term $2\epsilon \re ( |z_2|^2 z_3^2 \bar z_3^3  |z_4|^2)$ vanishes, and thus $ p$ is a sum of squares, making it automatically plurisubharmonic. Therefore, without the loss of generality, we can assume simultaneously that $z_2\neq0,$ $z_3\neq0,$ and $z_4\neq0.$ As a result, we can compute the Levi form in terms of vectors fields $$X=\sum_{j=2}^4 a_j \, z_j \, \frac{\partial}{\partial z_j}$$ and $\bar X,$ which keep the weight of each term of $p$ unchanged.

Rather than writing out the Levi form in full in terms of $X,\bar X$ as a quadratic form in $a_j,$ we observe that by Cauchy-Schwarz,
\begin{equation}
\label{cs1}
\left|2\re ( |z_2|^2 z_3^2 \bar z_3^3  |z_4|^2)\right| \le ( |z_2|^2 |z_3|^6 + |z_2|^2|z_3|^4 |z_4|^4)
\end{equation}
and 
\begin{equation}
\label{cs2}
\left|2\re ( |z_2|^2 z_3^2 \bar z_3^3  |z_4|^2)\right|  \le ( |z_2|^4 |z_3|^2 |z_4|^2 + |z_3|^8|z_4|^2).
\end{equation}

For the right-hand side expression in \eqref{cs1}, its kernel is given by the simultaneous vanishing of $\{a_2+3a_3 =0\}$ and $\{a_2+2a_3+2a_4 =0\}.$ For the right-hand side expression in \eqref{cs2}, its kernel is given by the simultaneous vanishing of $\{2a_2+a_3+a_4 =0\}$ and $\{4a_3+a_4 =0\}.$ The intersection of these kernels is just the origin, so for small $\epsilon,$ $p$ is indeed plurisubharmonic. Therefore, $r_0$ defines a pseudoconvex domain.



\bibliographystyle{alpha}


\begin{thebibliography}{KoMZ14}

\bibitem[BG77]{bloomgraham}
T. Bloom and I. Graham.
\newblock On ``type'' conditions for generic real submanifolds of {${\bf
  C}^{n}$}.
\newblock {\em Invent. Math.}, 40(3):217--243, 1977.

\bibitem[BS92]{bs}
H.~P. Boas and E.~J. Straube.
On equality of line type and variety type of real hypersurfaces in $\C^n$. 
{\em J. Geom. Anal.} {\bf 2} (1992), no. 2, 95--98. 

\bibitem[C84a]{C84}  
D.~W. Catlin.
Global regularity of the $\bar\d$-Neumann problem. Complex analysis of several variables (Madison, Wis., 1982), 39--49, Proc. Sympos. Pure Math., 41, Amer. Math. Soc., Providence, RI, 1984.

\bibitem[C84b]{catlinbdry}
D.~W. Catlin.
\newblock Boundary invariants of pseudoconvex domains.
\newblock {\em Ann. of Math. (2)}, 120(3):529--586, 1984.


\bibitem[C87]{C87} 
D.~W. Catlin.
Subelliptic estimates for the $\bar\d$-Neumann problem on pseudoconvex domains. {\em Ann. of Math. (2)}, 
126 (1): 131--191, (1987).


\bibitem[CD10]{fribourgcatda}
D.~W. Catlin and J.~P. D'Angelo.
\newblock Subelliptic estimates.
\newblock In {\em Complex analysis}, Trends Math., pages 75--94.
  Birkh\"auser/Springer Basel AG, Basel, 2010.

\bibitem[ChM74]{CM}
S.~S. Chern and J.~K. Moser.
\newblock Real hypersurfaces in complex manifolds.
\newblock {\em Acta Math.}, 133:219--271, 1974.

\bibitem[D82]{opendangelo}
J.~P. D'Angelo.
\newblock Real hypersurfaces, orders of contact, and applications.
\newblock {\em Ann. of Math. (2)}, 115(3):615--637, 1982.

\bibitem[Ko10]{martinIMRN}
M. Kolar.
The Catlin multitype and biholomorphic equivalence of models. 
{\em Int. Math. Res. Not. IMRN} 2010, no. 18, 3530--3548. 

\bibitem[KoM11]{martinfrancineCM}
M. Kolar and F. Meylan.
\newblock Chern-{M}oser operators and weighted jet determination problems.
\newblock In {\em Geometric analysis of several complex variables and related
  topics}, volume 550 of {\em Contemp. Math.}, pages 75--88. Amer. Math. Soc.,
  Providence, RI, 2011.


\bibitem[KoMZ14]{dimafrancinemartin}
M. Kolar, F. Meylan, and D. Zaitsev.
\newblock Chern-{M}oser operators and polynomial models in {CR} geometry.
\newblock {\em Adv. Math.}, 263:321--356, 2014.

\bibitem[K79]{kohnacta}
J.~J. Kohn.
\newblock Subellipticity of the {$\bar \partial $}-{N}eumann problem on
  pseudo-convex domains: sufficient conditions.
\newblock {\em Acta Math.}, 142(1-2):79--122, 1979.

\bibitem[KiZ17]{kimzaitsev}
S.Y. Kim and D. Zaitsev.
\newblock Jet vanishing orders and effectivity of {K}ohn's algorithm in
  dimension $3$. Preprint 2017.
\newblock https://arxiv.org/abs/1702.06908.

\bibitem[M92]{mcneal} 
J.~D.~McNeal.
Convex domains of finite type.
{\em J. Funct. Anal.} 108 (1992), no. 2, 361--373. 

\bibitem[N14]{ra}
A.~C. Nicoara.
\newblock {D}irect {P}roof of {T}ermination of the {K}ohn {A}lgorithm in the
  {R}eal-{A}nalytic {C}ase.
\newblock https://arxiv.org/abs/1409.0963v1.

\bibitem[S10]{siunote}
Y.-T. Siu.
\newblock Effective termination of {K}ohn's algorithm for subelliptic
  multipliers.
\newblock {\em Pure Appl. Math. Q.}, 6(4, Special Issue: In honor of Joseph J. Kohn. Part 2):1169--1241, 2010.

\bibitem[S17]{siunew}
Y.-T. Siu.
\newblock New procedure to generate multipliers in complex {N}eumann problem
  and effective {K}ohn algorithm.
\newblock {\em Sci. China Math.}, 60(6):1101--1128, 2017.

\bibitem[Y92]{jiyeyu}
J. Yu.
\newblock Multitypes of convex domains.
\newblock {\em Indiana Univ. Math. J.}, 41(3):837--849, 1992.

\bibitem[Z17]{dimatensors}
D. Zaitsev.
\newblock A geometric approach to {C}atlin's boundary systems.
\newblock Preprint 2017. 
https://arxiv.org/abs/1704.01808.

\end{thebibliography}


\end{document}